\documentclass{article}
\usepackage[left=3cm, right=3cm, top=2cm]{geometry}

\title{\bf Inverse Problems in Topological Persistence}
\author{{\bf Steve Oudot}  \\
  DataShape team, Inria\\
  \url{steve.oudot@inria.fr}
	\and 
	{\bf Elchanan Solomon} \\
	Department of Mathematics, Brown University \\
        \url{yitzchak_solomon@brown.edu}
}
\date{}

\usepackage{amsthm}
\usepackage[alphabetic]{amsrefs}

\usepackage{hyperref,url,amssymb,colonequals,pdfsync,float,mathtools,appendix,float,amsmath,tikz-cd} 

\begin{document}
	
	\maketitle


\newtheorem{theorem}{Theorem}[section]
\newtheorem{claim}[theorem]{Claim}
\newtheorem{conj}[theorem]{Conjecture}
\newtheorem{lemma}[theorem]{Lemma}
\newtheorem{example}{Example}
\newtheorem{prop}[theorem]{Proposition}
\newtheorem{problem}{Problem}
\newtheorem{corollary}[theorem]{Corollary}
\newtheorem{remark}{Remark}
\newtheorem{observation}[theorem]{Observation}
\newtheorem{counterexample}[theorem]{Counterexample}
\newtheorem{dictionary}[theorem]{Dictionary}
\newtheorem{definition}[theorem]{Definition}
\newcommand{\R}{\mathbb{R}}
\newcommand{\Z}{\mathbb{Z}}
\newcommand{\N}{\mathbb{N}}
\newcommand{\X}{\mathbb{X}}
\newcommand{\rtop}{\R\mathbf{Top}}
\newcommand{\kmod}{\mathbf{k}\text{-}\mathbf{Mod}}
\newcommand{\fieldk}{\mathbf{k}}
\newcommand{\isaac}[1]{{\color{magenta} \textsf{$\spadesuit\spadesuit\spadesuit$ isaac: [#1]}}}
\newcommand{\steve}[1]{{\color{blue} \textsf{$\clubsuit\clubsuit\clubsuit$ steve: [#1]}}}

\newcommand{\blue}[1]{{\color{black} #1}}
\newcommand{\magenta}[1]{{\color{magenta} #1}}
\newcommand{\green}[1]{{\color{ForestGreen} #1}}
\newcommand{\norm}[1]{ \left\lVert#1 \right\rVert}
\flushbottom 

\maketitle 

\tableofcontents 


	\begin{abstract}
	In this survey, we review the literature on inverse problems in topological persistence theory. The first half of the survey is concerned with the question of surjectivity, i.e. the existence of right inverses, and the second half focuses on injectivity, i.e. left inverses. Throughout, we highlight the tools and theorems that underlie these advances, and direct the reader's attention to open problems, both theoretical and applied.
\end{abstract}

\section*{Introduction} 

\addcontentsline{toc}{section}{Introduction} 

In recent decades, success in machine learning has revolved around the study of non-linear feature extraction and non-linear models. This paradigm uses large training sets and increased processing power to produce highly flexible models with ever increasing prediction accuracies. However, there is an emerging awareness among machine learning researchers and end-users that these non-linear techniques can be very hard to interpret. Often, the mapping from the input (data) space to the target (modeling) space is so complex that it is virtually impossible to predict what simple transformations in the target space might mean for real-world data, if they can be given any interpretation at all. Similarly, it is possible for slightly different input data sets to produce wildly divergent models. As prediction accuracy is only one part of the data analysis pipeline, many researchers are now studying the hard mathematical problems underlying the \emph{explainability} and \emph{interpretability} of machine learning algorithms.\\

The focus of this article is on Topological Data Analysis (TDA), which provides a set of feature extraction and modeling algorithms built around ideas and techniques from algebraic topology and metric geometry, and is particularly well-suited to studying data sets of complex shapes. Because of its origins in abstract mathematics, it is a prime candidate for modern research in explainability. In the following sections, we survey the work done in the TDA community on two topics of considerable interest: the \emph{preimage problem}, and \emph{discriminativity}.\\

\begin{figure}[htb]
	\centering
	\includegraphics[scale = 0.5]{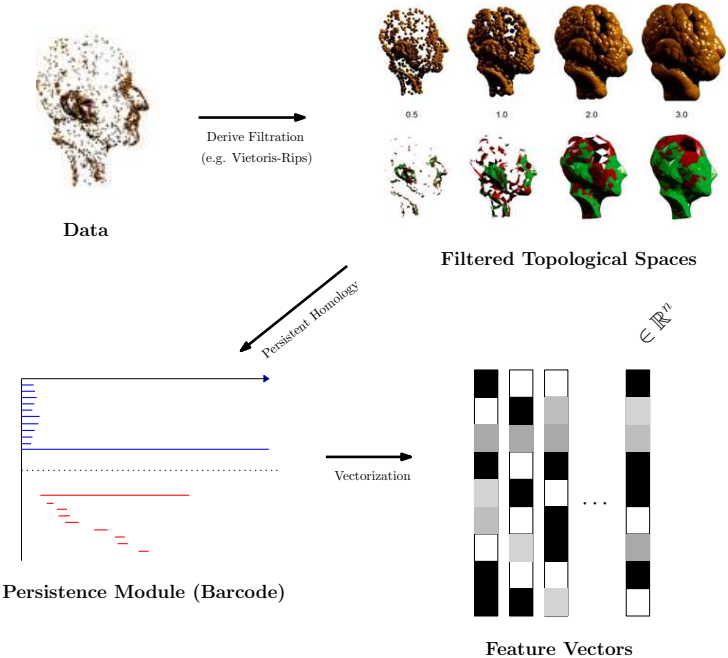}
	\caption{We are initially given a data set, such as a point cloud. From this, we derive a filtered topological space, using the Vietoris-Rips complex, \v{C}ech complex, $\alpha$-complex, etc. We then apply the homology functor to these topological spaces, obtaining a persistence module that can be represented as a barcode. Finally, for machine learning applications, there are methods of turning barcodes into feature vectors. Part of this figure has been adapted from~\cite[Fig. 6]{giesen2006conformal}.}
	\label{fig:pipeline}
\end{figure}

The central invariant of TDA is \emph{persistent homology}, which maps an input shape to a descriptor consisting of a set of intervals on the real line (called a \emph{barcode}). The persistent homology pipeline is outlined in Figure \ref{fig:pipeline}, which is explained in further detail in the background section.\\

 As it turns out, persistent homology provides provably stable descriptors, i.e. similar shapes are always mapped to similar sets of intervals. Insofar as explainability is concerned, it is natural to ask if every such descriptor corresponds to an input shape, and if so, how that shape might be approximated (the preimage problem). This is the focus of Section \ref{sec:right}, which considers positive preimage results for a variety of data types: model data, point clouds, and function-valued data. Another important question is that of injectivity: whether it is possible for two distinct shapes to produce identical descriptors, so that the resulting feature vectors cannot be used to distinguish between them. In Section \ref{sec:left}, we outline what is known about injectivity in the context of persistent homology, and consider two enriched feature-extraction models for which positive injectivity results have been proven.\\

As the reader may notice in the course of the following survey, explainability is not an accidental feature of the TDA framework, but an essential component of its philosophy. However, because TDA must often interface with traditional machine learning algorithms in practice, this survey is not independent of its analogues in the study of kernel machines and neural networks. Conversely, there is evidence that techniques of TDA can be used to study interpretability in other areas of data science, see \cite{chazal2013persistence,lee2017quantifying,gunnaronline} for examples. Thus, there is a promising two-way dialogue between researchers studying explainability in traditional machine learning and TDA.\\

{\bf Note}: As this is still an emerging topic in TDA, we need clarify that this survey is targeted at researchers working in applied topology and computational geometry. Otherwise, interested readers should have a working knowledge of elementary algebraic topology and homology theory (consult the appropriate chapters in~\cite{hatcher2005algebraic,munkres2018elements}  for a good introduction), as well as the essentials of commutative algebra (chapter 2 of Atiyah and McDonald's book \cite{atiyah1969commutative} is an excellent reference). For an introduction to the themes and tools of topological data analysis, the reader can consult the articles of Ghrist \cite{ghrist2008barcodes} and Carlsson \cite{carlsson2009topology}. Lastly, a more formal and comprehensive treatment can be found in the texts by Edelsbrunner and Harer \cite{edelsbrunner2010computational}, Ghrist \cite{ghrist2014elementary}, and Oudot \cite{oudot2015persistence}. 


\section{Background}
We now introduce definitions and constructions necessary for the rest of the survey. For us, persistence will be a functor from the category $\rtop$ to the category $\kmod$.

\begin{definition}
We define the category of $\mathbb{R}$-filtered topological spaces $\rtop$ to be the functor category from the poset category $(\mathbb{R},\leq)$ of ordered real numbers to the category $\mathbf{Top}$, whose objects are topological spaces, and whose morphisms are continuous maps. Additionally, throughout this paper, we stipulate that these continuous maps be set inclusions (following the TDA literature). Morphisms of $\mathbb{R}$-filtered topological spaces are then natural transformations between such functors.
\end{definition}

Put concretely, an object $X$ of $\rtop$ is a family of topological spaces $X(r)$ indexed by $r \in \mathbb{R}$, with set inclusions $X(r \leq s): X(r) \hookrightarrow X(s)$ for all $r \leq s \in \mathbb{R}$. A morphism of $\mathbb{R}$-filtered topological spaces $X$ and $Y$ is a family $\psi$ of continuous maps, $\psi(r): X(r) \to Y(r)$, with $\psi(s) \mid_{X(r)} = \psi(r)$ for $r \leq s$. Equivalently, we assert that the following square commutes.

\begin{center}
	\begin{tikzcd}[sep = large]
	X(r) \arrow[hookrightarrow]{r}{X(r \leq s)} \arrow{d}{\psi(r)} & X(s) \arrow{d}{\psi(s)}\\
	Y(r) \arrow[hookrightarrow]{r}{Y(r \leq s)}  &Y(s)
	\end{tikzcd}
\end{center}

A rich source of examples of $\mathbb{R}$-filtered topological spaces stems from \emph{point clouds} (see Figure \ref{fig:pipeline}). There are various ways to obtain an object in $\rtop$ from a point cloud $X$, with one of the most common being the Vietoris-Rips complex.

\begin{definition}
  Let $X \subset \mathbb{R}^d$ be a point cloud. The Vietoris-Rips (VR) filtration $VR(X)$ is a filtration on the full simplex on the set $X$ (i.e. the simplex of dimension $|X|-1$). For $r \in \mathbb{R}$, the subspace $\left(VR(X)\right)(r)$ consists of those simplices of diameter $\leq r$. We will call the diameter of a simplex $\tau$ its \emph{appearance time} in this filtration.
\end{definition}

One can also obtain $\mathbb{R}$-filtered topological spaces by using real-valued functions.

\begin{definition}
	Let $X$ be a topological space, and $f: X \to \mathbb{R}$ a continuous, real-valued function. We will write $(X,f)$ to denote the $\mathbb{R}$-filtered topological space consisting of the sublevel sets of $f$,
	\[(X,f)(r) = \{x \in X \mid f(x) \leq r\}.\]
\end{definition}

\begin{definition}
 We now define the category of \emph{persistence modules} $\kmod$ to be the functor category from the poset category $(\mathbb{R},\leq)$ to the category $\mathbf{Vect}$ of vector spaces over a fixed field $\fieldk$. Morphisms of persistence modules are then natural transformations between such functors.
\end{definition}
 
Put concretely, an object $M$ of $\kmod$ is a family of vector spaces $M(r)$ indexed by $r \in \mathbb{R}$, together with linear maps $M(r \leq s): M(r) \to M(s)$ for all $r \leq s \in \mathbb{R}$. These linear maps are required to satisfy the following compatibility axioms: $M(r \leq r) = \operatorname{id}_{M(r)}$, and $M(r \leq t) =  M(s \leq t) \circ M(r \leq s)$ for $r \leq s \leq t \in \mathbb{R}$. A morphism $\psi$ of persistence-modules $M$ and $N$ is a family of maps $\psi(r): M(r) \to N(r)$ making the following square commute for all $r \leq s$.
\begin{center}
\begin{tikzcd}[sep = large]
M(r) \arrow{r}{M(r \leq s)} \arrow{d}{\psi(r)} & M(s) \arrow{d}{\psi(s)}\\
N(r) \arrow{r}{N(r \leq s)}  &N(s)
\end{tikzcd}
\end{center}

We define the persistence map as follows.

\begin{definition}
Let $X$ be an $\mathbb{R}$-filtered topological space. The associated degree-$d$ persistence module $M$ has the degree-$d$ singular homology group $M(r) = H_{d}(X(r);\fieldk)$ at each index $r \in \mathbb{R}$, and the morphism $M(r \leq s): M(r) \to M(s)$  induced in homology by the inclusion $X(r) \hookrightarrow X(s)$ for each $r\leq s\in\R$. We will use the notation $PH_{d}(X) = M$ to indicate that $M$ is the degree-$d$ \emph{persistent homology} of $X$. When our $\mathbb{R}$-filtered topological space is the sublevel set filtration induced by a continuous real-valued function on a topological space, $f: T \to \mathbb{R}$, we will write $PH_{d}(T,f)$ for the resulting persistence module; this is called	 \emph{functional persistence} in the literature\footnote{Throughout the survey, we will use capital letters such as $X$ and $Y$ to refer to elements of both $\rtop$ and $\mathbf{Top}$. It will always be made clear, either explicitly or from the context, which one is intended.}.
\end{definition}

 For the remainder of the survey, we will omit any reference to the choice of field $\fieldk$, except when it is necessary to be explicit.

 \bigskip

 When computing homology in multiple degrees, we will want to keep track of all the resulting persistence modules at once. The appropriate algebraic object is a \emph{graded persistence module}.
\begin{definition}
 A graded persistence module $M = \bigoplus_{i \in \mathbb{N}} M_{i}$ is the direct sum of a family of persistence modules indexed over the natural numbers, together with the labeling that records which factor is associated to which number\footnote{Note that the grading here happens in the category of abelian groups, rather than in the category of modules. That is, the grading does not come with a multiplicative structure.}. The graded persistence module associated to an $\mathbb{R}$-filtered topological space $X$ is then 
	\[PH(X) = \bigoplus_{i \in \mathbb{N}} PH_{i}(X) \]
\end{definition}
    
    Though persistence modules are not vectors, they still live in a metric space. Indeed, the category $\kmod$ comes equipped with an extended pseudo-metric: the \emph{interleaving distance} $d_{I}$.
    
\begin{definition}
 An $\epsilon$-interleaving of persistence modules $M$ and $N$ consists of two families of morphisms, $f(r): M(r) \to N(r+\epsilon)$ and $g(r):N(r) \to M(r+\epsilon)$, making the following four diagrams commute for all $r \leq s$.
\begin{center}
    \begin{tikzcd}
    M(r) \arrow{r}{M(r \leq s)} \arrow{rd}{f(r)} & M_{s} \arrow{rd}{f(s)} &\\
     & N(r+\epsilon) \arrow{r}[below, yshift = -1ex]{N(r + \epsilon \leq s + \epsilon)} & N(s+\epsilon)
    \end{tikzcd}
    \begin{tikzcd}
    N(r) \arrow{r}{N(r \leq s)} \arrow{rd}{g(r)} & N(s) \arrow{rd}{g(s)} &\\
     & M(r+\epsilon) \arrow{r}[below, yshift = -1ex]{M(r + \epsilon \leq s + \epsilon)} & M(s+\epsilon)
    \end{tikzcd}    
\end{center}
\begin{center}
    \begin{tikzcd}
    M(r) \arrow{rd}{f(r)} \arrow{rr}{M(r \leq r+2\epsilon)} && M(r + 2\epsilon)\\
    & N(r+\epsilon) \arrow{ru}{g(r + \epsilon)} &
    \end{tikzcd}
    \begin{tikzcd}
    N(r) \arrow{rd}{g(r)} \arrow{rr}{N(r \leq r+2\epsilon)} && N(r + 2\epsilon)\\
    & M(r+\epsilon) \arrow{ru}{f(r + \epsilon)} &
    \end{tikzcd}    
\end{center}
\end{definition}
  
  Intuitively, one can think of such an interleaving as an approximate isomorphism of persistence modules. Indeed, a $0$-interleaving is exactly an isomorphism.

  \begin{definition}
    The interleaving distance $d_{I}$ between $M$ and $N$ is the infimum of values $\epsilon$ for which an $\epsilon$-interleaving exists. It satisfies the triangle inequality but can be zero between non-isomorphic modules, or equal to infinity.
\end{definition}

The category of persistence modules is abelian, which, among other things, allows one to take direct sums of persistence modules, defined pointwise.

\begin{definition}
Let $M$ and $N$ be a pair of persistence modules. We define their direct sum $M \oplus N$ to be the persistence module with vector spaces $(M \oplus N)(r) = M(r) \oplus N(r)$ and maps $(M \oplus N)(r \leq s) = M(r \leq s) \oplus N(r \leq s)$ for any $r \leq s$.
\end{definition}

An \emph{indecomposable} persistence module is one that cannot be written as the sum of two nonzero persistence modules. Examples of such modules include the \emph{interval persistence modules} $\fieldk_I$, defined as follows. Given an interval $I \subset \mathbb{R}$, let $\fieldk_{I}$ be such that $\fieldk_{I}(r) = \fieldk$ for $r \in I$ and has rank zero otherwise, and that $\fieldk_{I}(r \leq s) = id_{\fieldk}$ for $r \leq s \in I$ and is the zero map otherwise.\\

The category $\kmod$ contains some wild objects that are difficult to work with. Thus, it is necessary to restrict our attention to a  class of well-behaved persistence modules which suffices for practical applications:

\begin{definition}
We say that a persistence module $M$ is \emph{pointwise finite-dimensional} (pfd) if each vector space $M(r)$ is finite dimensional.
\end{definition}
 The following theorem asserts that every pfd persistence module has a particularly simple decomposition into indecomposables, and highlights the important role played by interval modules in the theory of persistence.

\begin{theorem}[\cite{crawley2015decomposition}]
\label{thm:crawleyboevey}
Every pfd persistence module is isomorphic to the direct sum of interval modules. Moreover, the decomposition is unique up to isomorphism and reordering of the terms.
\end{theorem}

From Theorem \ref{thm:crawleyboevey}, we see that pfd persistence modules admit a complete invariant: the \emph{barcode} formed by the collection of intervals involved in the direct sum decomposition of the module. More generally, we call a \emph{barcode} any multi-set of intervals. This terminology comes from plotting the intervals along a common axis, as in Figure \ref{fig:pipeline}.  

The space of barcodes has a natural metric: the \emph{bottleneck distance} $d_{B}$, defined as follows.

\begin{definition}
An $\epsilon$-matching between multi-sets of intervals $\mathcal{I}$ and $\mathcal{J}$ is a bijection between subsets $\mathcal{I}' \subseteq \mathcal{I}$ and $\mathcal{J}' \subseteq \mathcal{J}$ such that if the interval $[a,b] = I \in \mathcal{I}'$ is matched with the interval $[c,d] = J \in \mathcal{J}'$ then $\max \{|a-c|,|b-d|\} \leq \epsilon$, and such that any interval in $\mathcal{I} \setminus \mathcal{I}'$ or $\mathcal{J} \setminus \mathcal{J}'$ has diameter at most $2\epsilon$. The bottleneck distance between barcodes is the infimum of values $\epsilon$ for which there exists an $\epsilon$-matching between them.
\end{definition}

Persistent homology enjoys a variety of stability theorems. We recall here three of the most fundamental ones:

\begin{theorem}[Algebraic Stability, \cite{chazal2016structure,bauer2014induced,chazal2009proximity}]
For a pair $M$ of $N$ of pfd persistence modules with barcodes $B(M),B(N)$, the interleaving distance bounds the bottleneck distance.
\[d_{B}(B(M),B(N)) \leq d_{I}(M,N)\]
\end{theorem}

In fact, the above inequality is an equality, a result known as the \emph{isometry theorem}, cf. \cite{lesnick2015theory,chazal2016structure}.

\begin{theorem}[Geometric Stability, \cite{chazal2014persistence}]
Let $X$ and $Y$ be totally bounded metric spaces whose VR complexes have degree-$i$ persistence modules $M$ and $N$ respectively. If we let $B(M)$ and $B(N)$ denote the respective barcodes of these persistence modules, and $d_{GH}(X,Y)$ denote the Gromov-Hausdorff distance between these spaces, then
\[d_{B}(B(M),B(N)) \leq 2d_{GH}(X,Y).\]
\end{theorem}

\begin{theorem}[Functional Stability, \cite{chazal2016structure,cohen2005stability}]
	\label{thm:funcstab}
Let $X$ be a topological space, and let $f,g:X \to \mathbb{R}$ be two functions whose sublevel sets have finite-dimensional homology groups. Then $(X,f)$ and $(X,g)$ give rise to pfd functional persistence modules $M$ and $N$ with
\[d_{B}(B(M),B(N)) \leq \|f-g\|_{\infty}.\]
\end{theorem}

In the remainder of the survey, we will slightly abuse notation and write $d_{B}(M,N)$ in place of $d_{B}(B(M),B(N))$. 

\section{Persistence and Right Inverses}
\label{sec:right}

\subsection{Persistent Moore Spaces}

\begin{figure}[htb]
	\centering
	\includegraphics[scale = 0.55]{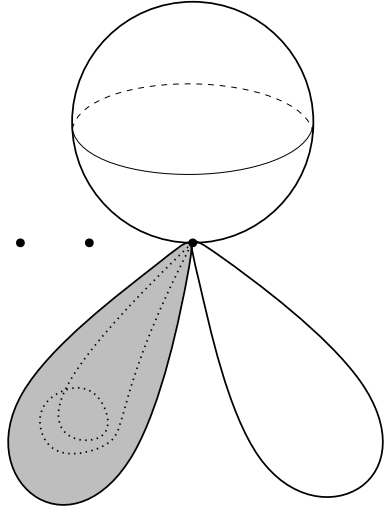}
	\caption{The space $X$ above is built using wedge sums and disjoint unions of Moore spaces. By construction, we have $H_{0}(X;\Z) \cong \Z^3$, $H_{1}(X;\Z) \cong \mathbb{Z} \oplus \mathbb{Z}/2\mathbb{Z}$, $H_{2}(X;\Z) \cong \mathbb{Z}$, and $H_i(X;\Z)=0$ for all $i\geq 3$.}
	\label{fig:moore}
\end{figure}

In this section, the input shapes of interest are $\mathbb{R}$-filtered topological spaces. Before addressing the existence of right inverses for persistent homology, let us review what is known for the usual (non-persistent) homology functor. It is a standard fact that homology admits a right inverse, in that every finitely-generated abelian group arises as the singular homology of some topological space $X$:
\begin{theorem}\label{thm:Moore}
	For any finitely generated graded abelian group $G=\bigoplus_{i \in \mathbb{N}} G_{i}$ such that $G_0$ is free and nontrivial, there is a topological space~$X$ such that $H_{i}(X;\Z) \cong G_{i}$ for all $i \in \mathbb{N}$.
\end{theorem}
%
Let us review the proof of this classical result, which proceeds by
constructing topological spaces realizing increasingly
varied groups in each degree $i>0$ separately. These spaces are called
{\em Moore spaces}, and we refer the reader to Section~2.2 (in
particular Example~2.40) in~\cite{hatcher2005algebraic} for a
background discussion. See also Figure~\ref{fig:moore} in this paper
for an illustration of the construction.\\

For the trivial group $G_i=0$, an appropriate Moore space is the one-point space $X_i = \{\ast\}$, which has trivial homology in all positive degrees. For the infinite cyclic group $G_i=\Z$, we can take $X_i =
\mathbb{S}^{i}$, the $i$-dimensional sphere. For a finite cyclic group
$G_i=\Z/n\Z$, we can glue the boundary of the disc $\mathbb{D}^{i+1}$
to the sphere $\mathbb{S}^{i}$ by a map of degree~$n$. In either case
we get a space~$X_i$ with degree-$i$ homology isomorphic to $G_i$ and
with trivial homology in the other positive degrees.\\

To realize an arbitrary finitely generated abelian group~$G_i$, we rely on the fact that such a group decomposes as a (finite) direct sum of cyclic groups:
\begin{equation}\label{eq:abelian_dec}
G_i\cong \bigoplus_{j=1}^{n_i} G_{i,j}, \mbox{ where each $G_{i,j}$ is cyclic.}
\end{equation}
Additionally, we make use of the following connection  between direct sums of homology
groups and wedge sums of spaces\footnote{This
	connection holds provided that the basepoints are chosen in such a way that
	they form good pairs with their associated spaces, which is the case
	here since all our spaces are CW-complexes.}:
\begin{equation}\label{eq:wedge}
\forall k>0,\quad H_k(\bigvee_{\alpha \in A} X_{\alpha};\Z) \cong \bigoplus_{\alpha \in A} H_{k}(X_{\alpha};\Z).
\end{equation}
This gives us a way to realize our group~$G_i$: given its
decomposition~(\ref{eq:abelian_dec}) and a collection of Moore
spaces $X_{i,j}$ realizing the cyclic summands~$G_{i,j}$, we take as
our Moore space the wedge sum $X_i=\bigvee_{j=1}^{n_i}X_{i,j}$,
which by~(\ref{eq:wedge}) has degree-$i$ homology isomorphic to
$G_i$ and trivial homology in the other nonzero degrees.\\

Finally, coming back to our initial graded group~$G$, we work with all homology degrees~$i>0$ at once and take the wedge sum
\[
Y = \bigvee_{i=0}^\infty X_i = \bigvee_{i=0}^\infty\ \bigvee_{j=1}^{n_i} X_{i,j},
\]
which by~(\ref{eq:wedge}) again has degree-$i$ homology isomorphic
to $G_i$ for each $i>0$. Since the whole space~$Y$ is path-connected
by construction, its degree-$0$ homology is isomorphic to~$\Z$, so
to complete the proof of Theorem~\ref{thm:Moore} we take~$X$ to be
the disjoint union of~$Y$ with $r-1$ copies of the one-point space,
where $r>0$ is the rank of the free group~$G_0$.\\

\begin{figure}[htb]
	\centering
	\includegraphics[scale=0.75]{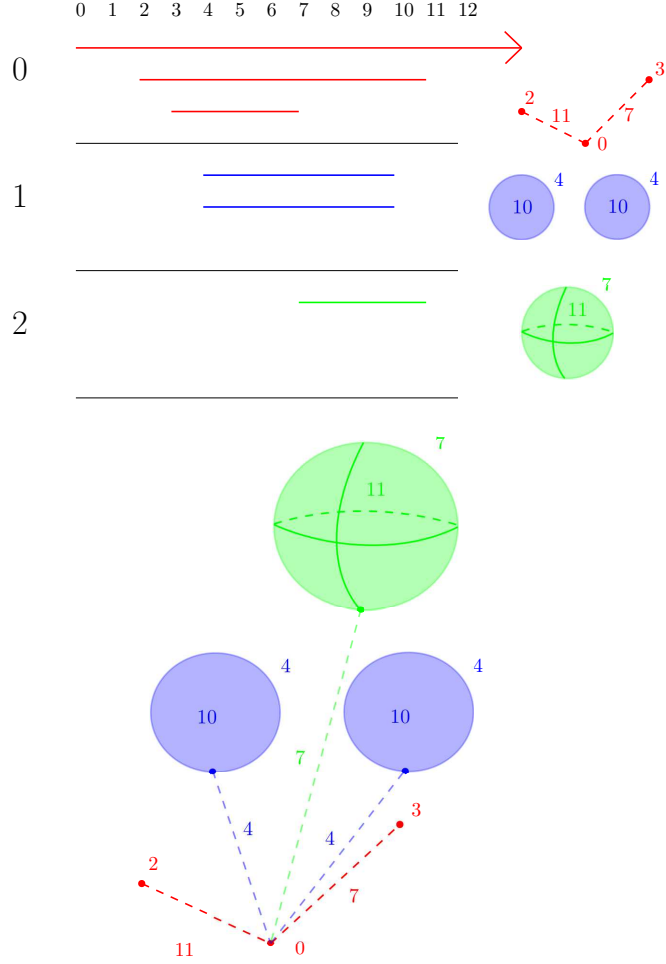}
	\caption{
		In the above figure, we construct persistent Moore spaces in each degree separately, and then glue them together to produce a right inverse for the entire graded persistence module.}
	\label{fig:persistentmooreconstruction}
\end{figure}

Transitioning to the case of (singular) persistent homology, we have
the following right-inverse theorem, where degree~$0$ again plays a
special role:
\begin{theorem}\label{thm:pers_Moore}
	Given a graded pfd persistence module $M=\bigoplus_{i \in
		\mathbb{N}} M_{i}$, if the barcode decomposition of~$M$ has a
	right-infinite interval in degree~$0$ that contains all the other
	intervals in the barcode (including all degrees), then there is an
	$\mathbb{R}$-filtered topological space $X$ with $PH(X) = M$.
\end{theorem}
Mirroring the construction in the non-persistent case, we begin by
realizing single interval modules in a single homology degree, then
we work our way up in complexity (see
Figure~\ref{fig:persistentmooreconstruction} for an illustration). To
that end, we introduce the concept of {\em persistent Moore space}
for an interval module:
\begin{definition}
	Given a homology degree $i>0$ and an interval $I \subseteq \mathbb{R}$, the
	\emph{persistent Moore space} $\mathbb{S}^{i}_{I}$
	is the following $\mathbb{R}$-filtered topological space, where the
	notation $r < I$ (resp. $r>I$) means that $r$ is less than
	(resp. greater than) every element~of~I:
	\[\mathbb{S}^{i}_{I}(r) =  \begin{cases} 
	\emptyset & r < I \\
	\mathbb{S}^{i} & r \in I \\
	\mathbb{D}^{i+1} & r > I
	\end{cases}\]
	Implicit in this formula is the fact that the boundary of the
	($i+1$)-disk is glued to the $i$-sphere by the identity map. 
\end{definition}
The persistent homology of $\mathbb{S}^{i}_{I}$ in degree $i$ is
isomorphic to the interval $I$-module, while it is trivial in the
other nonzero degrees. This construction thus produces an
$\mathbb{R}$-filtered topological space (persistent Moore space)
realizing any interval module in any fixed homology degree~$i>0$.\\

In order to extend the construction to arbitrary pfd persistence
modules, we use Theorem~\ref{thm:crawleyboevey} to decompose any such
module~$M$ into interval summands:
\[
M\cong\bigoplus_{j\in J} \fieldk_{I_j}.
\]
Then, given a fixed homology degree~$i>0$, for every interval
summand~$\fieldk_{I_j}$ we take a copy $X_{I_j}$ of the corresponding
persistent Moore space~$\mathbb{S}^{i}_{I_{j}}$. In order to combine
these spaces and realize the direct sum of their corresponding
interval modules, we choose a fixed basepoint on each copy of the
sphere~$\mathbb{S}^i$ and build a filtered version of the wedge sum, denoted by $\bigvee_{j\in J} X_{I_j}$ and defined as follows:
\begin{equation}\label{eq:filt_wedge}
\forall r\in\R,\quad (\bigvee_{j\in J} X_{I_j})(r) =  \bigvee_{j\in J} X_{I_j}(r),
\end{equation}
with the convention that $X\vee\emptyset = \emptyset\vee X = X$.  Note
that there are natural inclusions $X_{I_{j'}}(r)\hookrightarrow
\bigvee_{j\in J} X_{I_j}(r)$ for every $j'\in J$ and $r\in\R$,
so~(\ref{eq:filt_wedge}) yields a well-defined filtered
space. Moreover, it turns out that the isomorphism in~(\ref{eq:wedge})
is given by the direct sum of such inclusions, therefore, by
functoriality of homology, Eq.~(\ref{eq:wedge}) induces an isomorphism
between the degree-$i$ persistent homology of the filtered wedge
sum~(\ref{eq:filt_wedge}) and the persistence module $\bigoplus_{j\in J} \fieldk_{I_j} \cong M$.\\

As in the non-persistent setting, the graded version of this
construction works exactly the same way, by considering all homology
degrees~$i>0$ at once and taking the appropriate filtered wedge sum of
persistent Moore spaces. This yields an $\R$-filtered space~$Y$
whose graded persistent homology is isomorphic to a given graded pfd
persistence module~$M$, except possibly in degree~$0$.\\

Notice that the degree-$0$ persistent homology of~$Y$ is
isomorphic to a single interval module~$\fieldk_{I_0}$, since by
construction at each index~$r\in\R$ the space $Y_r$ is either empty or
path-connected. More precisely, $I_0$ is the smallest right-infinite
interval containing all the intervals in the barcode decomposition
of~$M$ in degrees $i>0$. Recalling now our assumption that the barcode
of~$M$ has a right-infinite interval $I'_0$ in degree~$0$ that
contains all these intervals (and therefore also~$I_0$), we want to
change the filtered space~$Y$ so that its degree-$0$ barcode now has a
single interval equal to~$I'_0$ while its degree-$i$ barcodes
remain unchanged for all $i>0$. This is done simply by taking the
filtered wedge sum of~$Y$ with the filtered one-point space
\[
P_r =  \begin{cases} 
\emptyset & r < I'_0 \\
\{\ast\} & r \in I'_0
\end{cases}\]

Finally, we can
further change~$Y$ into a filtered space~$X$ that has the
same barcode decomposition as~$Y$ in degrees~$i>0$ and that
acquires the missing intervals from the degree-$0$ barcode of~$M$. To do so, we take disjoint unions of $Y$ with filtered one-point spaces, giving rise to degree-$0$ bars with the appropriate left endpoints. We then specify, for each bar $I$, the corresponding right endpoint by gluing in a filtered edge that connects the one-point space associated with $I$ to the central connected component $Y$. We leave the details of this step as an exercise to the interested reader.

\begin{remark}
	 Our construction of persistent Moore spaces is a simplified version of the one introduced by Lesnick (\cite{lesnick2015theory}, \textsection~5.4); in exchange for dismissing the assumption that the spaces be compact, Lesnick's filtration arises from a real-valued function on a topological space. Moreover, Lesnick goes on to demonstrate a stronger result: for any pair of
	pfd persistence modules $M,N$
	and for any homology degree~$i$, there exists a common topological
	space $X$ and a pair of maps $\gamma^{M}, \gamma^{N}:X\to\R$, such
	that $PH_{i}(X,\gamma^{M})\cong M$ and $PH_{i}(X,\gamma^{N}) \cong
	N$. Lesnick further shows that the distance
	$d_{\infty}(\gamma^{M},\gamma^{N})$ between the maps can be made
	arbitrarily close to the distance $d_{I}(M,N)$ between the
	persistence modules.
\end{remark}

\subsection{Point Cloud Continuation}

We have seen that the map $PH:\rtop \to \kmod$ is surjective. However, data often takes the form of \emph{point clouds}: finite subsets of $\mathbb{R}^d$. The techniques of persistence homology can be applied to filtered spaces derived from these point clouds, such as their Vietoris-Rips (VR), \v{C}ech, or $\alpha$-filtrations. It is then natural to ask about the right-inverse problem for point clouds. Namely:
\begin{itemize}
    \item For $i \geq 0$, does every persistence module arise as the degree-$i$ VR/$\alpha$-complex persistence of a point cloud?
    \item If a given persistence module \emph{does} comes from a point cloud, can that point cloud be computed effectively?
\end{itemize}

The first question has a negative answer. To give a simple example, every zero-dimensional homology class of the VR filtration of a point cloud is born at zero, and hence any persistence module containing an interval summand born after zero cannot come from the zero-dimensional persistence of a point cloud. In general, it is unknown how to determine which persistence modules come from point clouds. However, the second question, that of computation, can sometimes be answered in the affirmative, at least locally, by using a continuation method.\\

\begin{figure}[htb]
    \centering
    \includegraphics[scale=0.8]{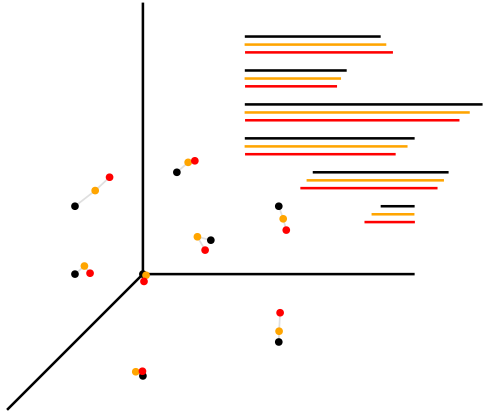}
    \caption{Modifying a persistence module via point cloud continuation.}
    \label{fig:pcc}
\end{figure}

The approach adopted in \cite{gameiro2016continuation} is the following. One is given an initial point cloud $P$ together with the the persistence module $M = PH(VR(P))$ induced in homology by its VR filtration. One then specifies a target persistence module $M'$ that is believed to be (close to) the persistence module of some unknown point cloud $P'$. The idea is then to use the Newton-Raphson method to make successive adjustments to $P$, incrementally bringing its persistence module closer to $M'$, as in Figure \ref{fig:pcc}.

To justify this approach, one must first make sense of the persistent homology algorithm as an actual function, with a well-defined domain and co-domain. To do this, let us consider the following segmentation of the VR persistence algorithm:

\begin{enumerate}
    \item A point cloud $X$ of $n$ points in $\mathbb{R}^d$, ordered as $x_1, \cdots, x_n$, can be associated with a vector in $\mathbb{R}^{nd}$.
    
    \item An ordering of the points $X$ gives rise to a lexicographc ordering on the powerset $2^X$, which allows us to associate a filtered simplicial complex on $X$ with a vector in $\mathbb{R}^{2n-1}$ (we ignore the empty simplex).
    
    \item There is a map $S: \mathbb{R}^{nd} \to \mathbb{R}^{2^n - 1}$ that sends an ordered point cloud to its associated filtered simplicial complex. The real value associated to each simplex is its appearance time in the VR filtration, which is half the distance between its furthest pair of vertices (their corresponding edge is called the \emph{attaching edge} of the simplex).

    \item There is a permutation $\pi$ on the set of simplices on $X$, depending on the ordering of their appearance times (and the lexicographic ordering, to break ties), that orders simplices via the pairing coming from persistent homology. That is, a simplex giving birth to a homological feature is followed by the simplex that kills that same feature, except for the first simplex, a vertex giving rise to the $0$-dimensional homology class with infinite persistence. There is a corresponding linear map $R_{\pi}:\mathbb{R}^{2^n - 1} \to \mathbb{R}^{2^n -1}$ that applies this permutation to the standard basis vectors.

    \item There is a projection map $P: \mathbb{R}^{2^n -1} \to \mathbb{R}^{m}$, which kills off all pairs with zero persistence. This corresponds to a barcode with $k$ bars, where $m = 2k-1$ (the infinite bar corresponds to a single simplex). As with $R_{\pi}$, the map $P$ depends on the point cloud $X$.
    
    \item Taken all together, on the fixed point cloud $X$ of $n$ points in $\mathbb{R}^d$, the persistence map agrees with the map $P \circ R_{\pi} \circ S: \mathbb{R}^{nd} \to \mathbb{R}^{m}$.
    
\end{enumerate}

To show differentiability of the persistence map at $X$, the key observation is that it agrees with the map $P \circ R_{\pi} \circ S$ in a \emph{neighborhood} of $X$ (they will certainly not agree on all of $\mathbb{R}^{nd}$). However, there is a caveat: if the pairwise distances in $X$ are not all distinct, then, for another point cloud $X'$ arbitrarily close to $X$, it is possible that the pairing of critical simplices may be different, and indeed the appropriate permutation and projection maps may be different from $R_{\pi}$ and $P$. To address this problem, Gameiro et al. assume that the point cloud $X$ is in VR-general position:
\begin{itemize}
    \item Condition A: All of the points in $X$ are distinct.
    \item Condition B: All of the appearance times of edges are distinct. Equivalently, all the pairwise distances between points are distinct.
\end{itemize}

Condition B ensures that the ordering on simplices coming from their appearance times is stable in a small neighborhood of $X$ in $\mathbb{R}^{nd}$, as appearance times of edges are continuous functions of distances between the points. This means that the permutation $\pi$ for $X$ in the above pipeline will give the correct pairing for nearby point clouds $X'$; similarly, the projection $P$ for $X$ will also drop the zero persistence pairs for $X'$. Thus, in a small neighborhood of $X$, the persistence map agrees with the same map $P \circ R_{\pi} \circ S$.
As linear maps, $R_{\pi}$ and $P$ are clearly $C^{\infty}$ differentiable. The chain rule tells us that $P \circ R \circ S$ will be $C^{\infty}$ differentiable if $S$ is. For a pair of points $x,y$ in $\mathbb{R}^d$, $S$ assigns their corresponding edge $E$ an appearance time of $r(x,y) = \frac{1}{2} \norm{x-y}$. This has partial derivatives
\[\frac{\partial r}{\partial x} = \frac{1}{2} \frac{x-y}{\norm{x-y}}, \,\,\,\,\,\,\,\,\,\,\,\,\,\,\,\, \frac{\partial r}{\partial y} = \frac{1}{2} \frac{y-x}{\norm{x-y}}\]

Condition A above guarantees that this derivative is defined, and indeed that $r(x,y)$ is $C^{\infty}$. Moreover, since every simplex in our VR filtration appears with a certain attaching edge, and since condition B ensures that this attaching edge remains the same for nearby point clouds, we know that all components of the map $S$ are $C^{\infty}$. Thus, since the persistence map agrees with $P \circ R_{\pi} \circ S$ in a neighborhood of $X$, it, too, is $C^{\infty}$. However, this is not sufficient for the implementation of the standard Newton-Raphson method, which requires the Jacobian of the map to be invertible. To cope with this, the authors use a slightly modified iteration scheme based on the (Moore-Penrose) pseudo-inverse of the Jacobian. We remind the reader that the psuedo-inverse $A^\dagger$ of a matrix $A$ is characterized by the following axioms:
\[AA^{\dagger}A = A\]
\[A^{\dagger}AA^{\dagger} = A^{\dagger}\]
\[(AA^{\dagger})^{T} = AA^{\dagger}\]
\[(A^{\dagger}A)^{T} = A^{\dagger}A\]

If $A$ has SVD decomposition $A = V\Sigma W^T$ then we have $A^\dagger = W \Sigma^{\dagger} V^T$, where $\Sigma^\dagger$ is obtained from $\Sigma$ by inverting the nonzero diagonal elements.

\begin{theorem}[\cite{gameiro2016continuation}, Corollary of Proposition 4.2]
	\label{thm:gameiroconvergence}
When $m=nd$, the iteration scheme described above converges to a point cloud $X'$. Moreover, when the Jacobian of the persistence map has full rank at $X'$, we may conclude that $PH(VR(X')) = M'$, the target persistence module.
\end{theorem}

\subsection*{Functional Optimization and Continuation}

	In \cite{poulenard2018topological}, Poulenard et al. consider the following problem, similar to that studied in \cite{gameiro2016continuation}. One is given a simplicial complex $X$ and a real-valued function $f_{\alpha}: X \to \mathbb{R}$ which depends on a continuous parameter $\alpha$. The persistence module $M = PH(X,f_{\alpha})$ is then stored as a multi-set of intervals $\{(b_i,d_i)\}_{i}$. Finally, there is a real-valued functional $\mathcal{F}$ which takes this multi-set as input. For example, this functional might record the distance between $M$ and a target module $N$. Our goal is to optimize the functional $\mathcal{F}$ as a function of the parameter $\alpha$, with the challenge being that the pipeline from $\alpha$ to $\mathcal{F}$ incorporates the procedure of taking persistent homology. Thus, when applying the chain rule, it will be necessary to differentiate the endpoint values $b_{i}$ and $d_{i}$ with respect to $\alpha$, which is not clearly defined. Using an argument similar to, and often simpler than, that of \cite{gameiro2016continuation}, Poulenard et al. show how to locally associate such an endpoint value with a fixed vertex $v_{i} \in X$, so that the derivative $\partial b_{i}/\partial \alpha$ or $\partial d_{i}/\partial \alpha$ that shows up in the chain rule is locally replaced with $\partial f/ \partial \alpha \mid_{v_{i}}$. This allows one to locally define the gradient $\nabla_{\alpha} \mathcal{F}$, and so approximate an optimum via gradient descent.\\

	In the case of minimizing the distance between $M = PH(X,f_{\alpha})$ and a target persistence module $N$, Poulenard et al. do not provide any convergence guarantees analogous to Theorem \ref{thm:gameiroconvergence}. For other applications, they prove and make use of another inverse-type result in applied topology.
	
	\begin{definition}
		For a simplicial complex $X$, let $F(X)$ be the space of real-valued functions on $X$. For a pair of simplicial complexes $X$ and $Y$, a function $T: X \to Y$ induces a pullback linear transformation $T_{F}: F(Y) \to F(X)$ via precomposition.
	\end{definition}
	
	\begin{theorem}[Thm. 1 in \cite{poulenard2018topological}]
		\label{thm:point2point}
		 An invertible linear functional map $T_{F}: F(Y) \to F(X)$ corresponds to a
		\emph{continuous} bijective point-to-point map $T: X \to Y$ if and only if both $T_F$ and
		its inverse preserve pointwise products of pairs of functions, and
		moreover both $T_F$ and its inverse preserve the persistence modules
		of all real-valued functions. In other words:

		\[d_{B}(PH(Y,f),PH(X,T_{F}(f)) =0, \,\,\,\, \forall f.\]
		
		Note that preservation of products ensures that $T_F$ corresponds to
		a point-to-point map, whereas preservation of persistence modules
		guarantees that the underlying map is continuous.
	\end{theorem}

	Let us illustrate one application of Theorem \ref{thm:point2point}: improving continuity in functional maps. Poulenard et al. pick as functions $f_{i}$ the characteristic functions of certain connected components of $Y$, and define an energy on the space of linear functional maps $\{T_{F}: F(Y) \to F(X)\}$:
	
	\[E = \sum_{f_{i}} d_{B}(PH(Y,f_{i}),PH(X,T_{F}(f_{i}))\]
	
\begin{figure}[htb]
	\centering
	\includegraphics[scale = 0.6]{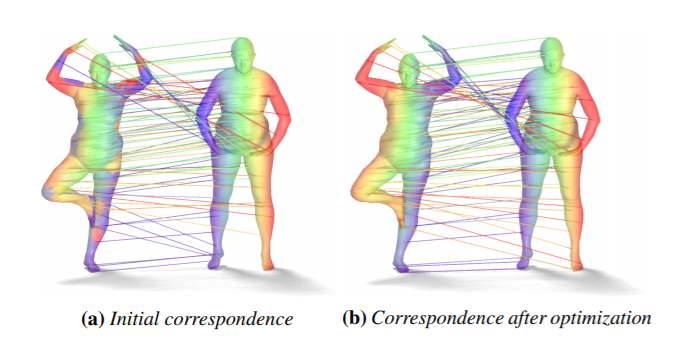}
	\caption{(a) A noisy functional map between two shapes converted to a point-to-point map. (b) The same map after topological optimization to improve continuity. The colors encode the $x$-coordinate function on the target shape and its pull-back on the source. Reproduced from \cite[\S 6]{poulenard2018topological}.}
	\label{fig:upgrade}
	
\end{figure}

 Theorem \ref{thm:point2point} implies that, if all connected components are taken in the above sum, and a zero-energy minimizer $T_{F}$ exists for the resulting functional, it corresponds to a continuous point-to-point map $T: X \to Y$. Optimizing this energy requires parametrizing the space of linear functional maps, which can be done using eigenfunctions of the Laplace-Beltrami operator to produce bases for $F(X)$ and $F(Y)$, and taking derivatives with respect to persistence modules, which they have already locally defined. As one can see in Figure \ref{fig:upgrade}, this can be used to upgrade poor initial correspondences between shapes to more continuous alignments.
	
\section{Persistence and Left Inverses}
\label{sec:left}

In what follows, the shapes of interest are metric spaces, and injectivity is considered \emph{up to isometry}. This presents a challenge for the existence of a left inverse to standard persistent homology invariants (functional, VR, $\alpha$-filtration, etc.), which are generally not sensitive enough to capture all this geometric data.
The following examples demonstrate some of the ways in which persistence maps can fail to be injective.
\begin{itemize}
    \item Rotating and translating a point cloud in $\mathbb{R}^d$ does not affect the persistent homology of its VR filtration (the same is true for the $\alpha$- or \v{C}ech filtration).
    \item The persistent homology of the VR or $\alpha$-filtration of a point cloud can also be preserved by non-isometries. Consider the three-point metric space $P_{\theta}$ obtained by taking the vertices of the triangle in Figure \ref{fig:threepoints} below. For any choice of $\theta \in [\pi/2,\pi]$, the persistence module of its VR filtration is the same (idem for the $\alpha$- or \v{C}ech filtration).
    \begin{figure}[htb]
        \centering
        \includegraphics[scale = 0.7]{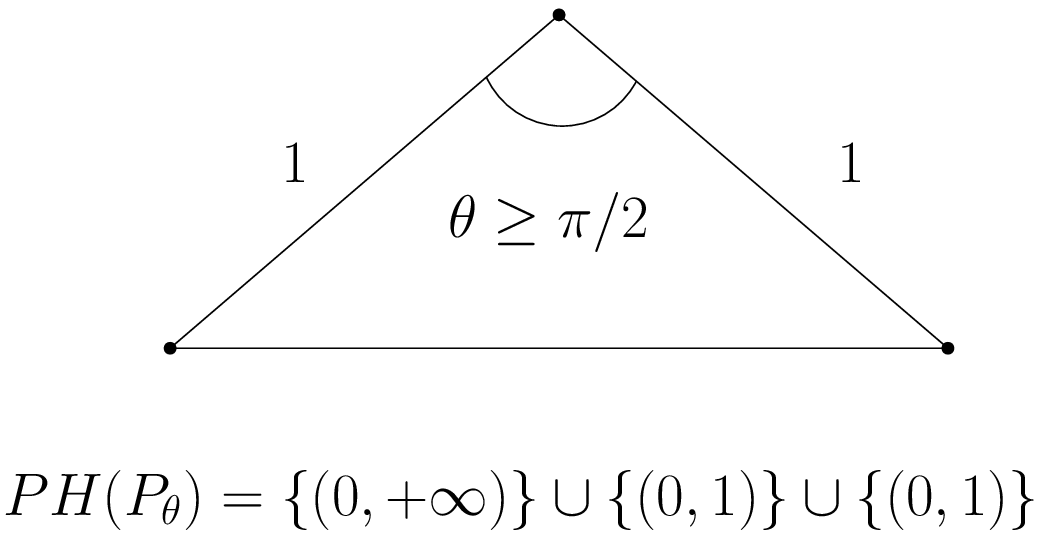}
        \caption{A family of non-isometric point clouds with the same persistence module.}
        \label{fig:threepoints}

\end{figure}
    \item Injectivity can also fail for intrinsic metric spaces. Indeed, the persistence module of the \v{C}ech filtration is identical for every geodesic tree, see e.g. Lemma~2 in \cite{gasparovic2018complete}. The same fact holds for the VR filtration, as shown in Appendix~A of \cite{BarcodeEmbeddings}.

    \item The following example gives insight into why the persistence
      map generally fails to be injective. Let $h: X \to Y$ be a
      homeomorphism of two, not-necessarily-isometric, topological
      spaces. Given a function $f: Y \to \mathbb{R}$, one can pull it
      back to a function $h^{\ast}(f) = f \circ h$ on $X$. By
      construction, $(X,h^{\ast}(f))$ and $(Y,f)$ have the same
      persistent homology. We see then that persistent homology is
      invariant under taking homeomorphisms, a much larger class of
      transformations than isometries.
    
    \item In \cite{curry2017fiber}, Curry characterized the fiber of the persistence map for functions on the unit interval, describing precisely which functions produce the same persistence module. However, in most settings, this is a hard, open problem.

\end{itemize}

These examples suggest that, to produce discriminative invariants using persistence, we must capture more information than a single persistence module.\\

\subsection{Extrinsic Persistent Homology Transforms}

\begin{figure}[htb]
    \centering
    \includegraphics[scale=0.6]{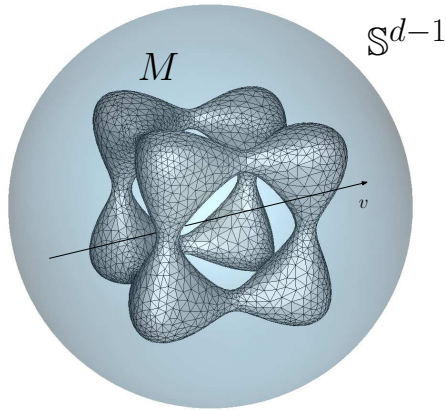}
    \caption{The map $f_{v}$.}
    \label{fig:PHT}
\end{figure}

In \cite{turner2014persistent}, Turner et al. propose the \emph{Persistent Homology Transform (PHT)}. The input to the PHT is a subanalytic compact subset $M$ of $\mathbb{R}^d$. For every direction $v \in \mathbb{S}^{d-1}$, one considers the function $f_{v}:M \to \mathbb{R}$, given by $f_{v}(x) = v \cdot x$ (see Figure \ref{fig:PHT}). The output of the PHT is then the map $PHT(M): \mathbb{S}^{d-1} \to \mathcal{M}$, sending a unit vector $v$ to the persistence module $PH(S,f_v)$.\\

Another, simplified invariant considered in \cite{turner2014persistent} is the \emph{Euler Characteristic transform (ECT)}. This is similar to the PHT, but instead of recording the sublevel-set persistence of the functions $f_{v}$, one computes their Euler Characteristic curves:
\[EC(f_v)(t) = \chi (\{x \in M \mid f_{v}(s) \leq t\})\]

If one writes $\{ \mathbb{R} \to \mathbb{Z}\}$ for the space of integer-valued functions on the real line, then this transform is the map $ECT:\mathbb{S}^{d-1} \to \{\mathbb{R} \to \mathbb{Z}\}$ that sends a unit vector $v$ to the function $EC(f_v)$. The object $ECT(M)$ lives in a Hilbert space, making it amenable to methods in classical statistics and machine learning. Indeed, Turner et al. show how to use the ECT to turn a set of meshes into a likelihood model on the space of embedded simplicial complexes. More precisely, they prove that, for $d = 2,3$, both of these transforms are injective, and hence provide  \emph{sufficient statistics} for probability measures on the space of linearly embedded simplicial complexes. Moreover, they provide an explicit algorithm to reconstruct $M$ from $PHT(M)$.\\

 Recent work of Ghrist et al.~\cite{ghrist2018persistent} and, independently, of Curry et al.~\cite{curry2018many}, using ideas of Schapira~\cite{schapira1995tomography}, demonstrates the injectivity of the ECT in all dimensions, and for the larger class of subanalytic compact sets. Because the Euler Characteristic curve of the functions $f_{v}$ can be derived from their persistence module, this, in turn, implies the injectivity of the PHT. These proofs of injectivity use the theory of constructible functions and Euler-Radon transforms, circumventing the involved, constructive arguments used in \cite{turner2014persistent}. Following \cite{curry2018many}, we introduce the necessary definitions and outline the proof below.\\
 
 Let $X$ be a real analytic manifold, and write $CF(X)$ for the space of \emph{constructible} functions on $X$. These are $\mathbb{Z}$-valued functions with subanalytic and locally finite level sets.

 \begin{definition}
 For a function $\phi \in CF(X)$, we define its Euler integral to be
     \[\int_{X} \phi(x) d\chi = \sum_{m \in \mathbb{Z}} m \, \chi(\{x \in X \mid \phi(x) = m\})\]
 \end{definition}

 \begin{definition}
 A morphism $f: X \to Y$ of real analytic  manifolds induces a \emph{pullback map} $f^{*}: CF(Y) \to CF(X)$ defined by $(f^{*}\phi)(x) = \phi(f(x))$ for $\phi \in CF(Y)$.
 \end{definition}
 \begin{definition}
 A morphism $f: X \to Y$ of real analytic manifolds induces a \emph{pushforward map} $f_{*}: CF(X) \to CF(Y)$ defined by $(f_{*}\phi)(y) = \int_{X} \phi\, 1_{f^{-1}(y)} \, d\chi$ for $\phi \in CF(X)$.
 \end{definition}
 
 These operations, taken together, allow us to define the following topological transform.
 
 \begin{definition}
 Let $S \subset X \times Y$ be a locally closed subanalytic subset of the product of two real analytic manifolds. Let $\pi_{X}$ and $\pi_{Y}$ be the projections from $X \times Y$ onto each of its factors. The \emph{Radon transform with respect to $S$} is the group homomorphism $R_{S}:CF(X) \to CF(Y)$ defined by $R_{S}(\phi) = (\pi_{Y})_{*}[ (\pi_{X})^{*}(\phi)1_{S}]$ for $\phi \in CF(X)$.
 \end{definition}
 
Schapira~\cite{schapira1995tomography} provides the following inversion theorem.

\begin{theorem}[Thm. 3.1 in \cite{schapira1995tomography}]
\label{thm:schapiroinversion}
Let $S \subset X \times Y$ and $S' \subset Y \times X$ define a pair of Radon transforms $R_{S}:CF(X) \to CF(Y)$ and $R_{S'}: CF(Y) \to CF(X)$. Denoting by $\overline{S}$ and $\overline{S'}$ the closure of these subsets, suppose that the projections $\pi_{Y}: \overline{S} \to Y$ and $\pi_{X}:\overline{S'} \to X $ are proper. Suppose further that there exists $\chi_1, \chi_2 \in \mathbb{Z}$ such that, for any $x \in X$, the fibers $S_{x} = \{y \in Y : (x,y) \in S\}$ and $S'_{x} = \{y \in Y: (y,x) \in S'\}$ satisfy the following criterion:
\[\chi(S_x \cap S'_{x}) =
\begin{cases}
	\chi_1 & \mbox{ if } x=x'\\
	\chi_2 & \mbox{ if } x\neq x'
\end{cases}
\]

Then for all $\phi \in CF(X)$,
\[(R_{S'} \circ R_{S})(\phi) = (\chi_1 - \chi_2)\, \phi + \chi_2 \left( \int_{X} \phi\, d\chi \right) 1_{X}\]
\end{theorem}

In particular, if $\chi_1 \neq \chi_2$ then the scaling term in $(R_{S'} \circ R_{s})$ is constant and nonzero. To take advantage of this theorem, \cite{ghrist2018persistent} define a Radon transform that can be computed using the ECT, and then find an appropriate ``inverse" Radon transform.\\

Let $X= \mathbb{R}^d$ and $Y = \operatorname{AffGr}_{d}$, the affine Grassmanian of hyperplanes in $\mathbb{R}^d$. Let $S \subset X \times Y$ be the set of pairs $(x,W)$, where the point $x$ sits on the hyperplane $W$. Letting $1_{M}$ be the indicator function of a bounded subanalytic subset $M \subset \mathbb{R}^d$, and $\pi_1$ and $\pi_2$ the projections of $X \times Y$ onto $X$ and $Y$ respectively, we compute:
\begin{align*}
    (R_{s}1_{M})(W) & = (\pi_{2})_{*}[(\pi_1^{*}1_M)1_S](W)\\
    &= \int_{(x,W) \in S} (\pi_{1}^{*}1_M)d\chi\\
    &= \int_{x \in M \cap W} d\chi\\
    &= \chi(M \cap W)
\end{align*}

To see that $\chi(M \cap W)$ can be computed from the ECT, let $W$ be defined by some unit vector $v$ and scalar $t$, i.e. $W = \{x : x \cdot v = t\}$. Then, using the inclusion-exclusion property of the Euler characteristic:
\begin{align*}
    \chi(M \cap W) &= \chi(\{x \in M: x \cdot v = t\})\\
     &= \chi(\{x \in M: x \cdot v \leq t\} \cap \{x \in M: x \cdot (-v) \leq -t\})\\
     &= \chi(\{x \in M: x \cdot v \leq t\}) + \chi(\{x \in M: x \cdot (-v) \leq -t\})\\ & - \chi(M)\\
     &= ECT(M)(v,t) + ECT(M)(-v,-t) - ECT(M)(v,\infty),
\end{align*}
 where $ECT(M)(v,\infty)$ is defined to be $\displaystyle\lim_{t \to +\infty} ECT(M)(v,t)$, which converges to $\chi(M)$ when $M$ is bounded.\\

Thus, if the Radon transform $R_{S}$ is injective, so is the ECT, as if $ECT(M) = ECT(M')$ for a pair of subanalytic subsets $M, M' \subset \mathbb{R}^d$ then $R_{S}1_{M} = R_{S}1_{M'}$. What remains to be shown, then, is that $R_{S}$ is indeed injective. We take $S' \subset Y \times X$ to consist of pairs $(W,x)$ where $x$ lies on the hyperplane $x'$. To apply Theorem \ref{thm:schapiroinversion}, we consider the intersection of fibers in $S$ and $S'$. For a fixed $x \in X$, $S_{x} = S'_{x} \subset Y$ is the set of hyperplanes passing through $x$, which is homeomorphic to the projective space $\mathbb{R}P^{d-1}$, which has Euler characteristic
\[\chi_1 = \chi(S_x \cap S^{'}_{x}) = \chi(\mathbb{R}P^{d-1}) = \frac{1}{2}(1 + (-1)^{d-1})\]

For a pair of distinct points $x \neq x'$, the intersection of fibers $S_{x} \cap S^{'}_{x'} \subset Y$ consists of all hyperplanes intersecting both of these points, a subset homeomorphic to $\mathbb{R}P^{d-2}$. Thus
\[\chi_2 = \chi(S_x \cap S^{'}_{x'}) = \chi(\mathbb{R}P^{d-2}) = \frac{1}{2}(1 + (-1)^{d-2})\]

By Theorem \ref{thm:schapiroinversion}, 
\[(R_{S'} \circ R_{S})(1_M) = (-1)^{d-1}1_{M} + \frac{1}{2}(1 + (-1)^{d-2})\chi(M)1_{\mathbb{R}^d}\]

Thus, if $R_{S}1_M = R_{S}1_{M'}$, then, composing with $R_{S'}$ and applying the above formula and rearranging terms, we obtain:
\[(-1)^{d-1}(1_{M} - 1_{M'}) = \frac{1}{2}(1 + (-1)^{d-2})(\chi(M')-\chi(M))1_{\mathbb{R}^d} \]

The right-hand side is a constant function, and so the left-hand side must be too. The difference of two non-zero indicator functions is constant precisely when it is equal to zero, so that $1_{M} = 1_{M'}$ and hence $M = M'$, demonstrating injectivity.

\subsection*{How many directions suffice?}

The injectivity results of \cite{turner2014persistent,ghrist2018persistent} require us to compute the PHT or ECT for \emph{every} vector on the sphere $\mathbb{S}^{d-1}$. Thus it is natural to ask if injectivity can be obtained with only finitely many directions. We should clarify that we are not asking for finitely many fixed directions to distinguish an infinite family of shapes. Rather, we would like to know if the identity of a given subanalytic set $S$ can be inferred by computing and comparing the PHT or ECT along a finite sequence of directions, with these directions being chosen in real time. There are two positive results in this vein, both restricted to the case of simplicial complexes, rather than arbitrary subanalytic sets.\\

The first result is that of \cite{belton2018learning}, specifically for the case of planar graphs. They demonstrate how to use three directions on the circle $\mathbb{S}^1$ to determine the location of the vertices of a planar graph $S$. The first two direction vectors are $(1,0)$ and $(0,1)$, and the third direction can be computed using the persistence modules derived from the first two. If $S$ has $n$ vertices, this vertex-localizing algorithm runs in $O(n \log n)$ time. Once the locations of the vertices are identified, one tests for the existence of an edge between pairs of vertices by using another three persistence modules (the directions of which are derived from the locations of the vertices). This pair-wise checking for edges introduces a quadratic term into the running time:
\begin{theorem}[Thm. 11 in \cite{belton2018learning}] Let $M$ be a linear plane graph with $n$ vertices. The vertices, edges, and exact embedding of $M$ can be determined using persistence modules along $O(n^2)$ different directions.
\end{theorem}

The second result, proved in \cite{curry2018many}, applies to finite, linearly embedded simplicial complexes $S \subset \mathbb{R}^d$ for any dimension~$d$. However, their bound on the number of directions is not simply a function of the number of vertices in $S$, but also of its geometry. In particular, it depends on the following three constants.
\begin{itemize}
    \item $d$ -- the embedding dimension.
    \item $\delta$ -- a constant with the following property: for any vertex $x \in M$ there is a ball $B$ of radius $\delta$ in the sphere $\mathbb{S}^{d-1}$, such that that for all $v \in B$ the Euler curve of $f_v$ changes values at $t = v \cdot x$. If one works with the PHT instead of the ECT, the analogous requirement is that the persistent homology coming from $f_v$ has an off-diagonal point with birth or death value $v \cdot x$. These conditions ensure that the vertex $x$ is \emph{observable} for the ECT or PHT in some simple set of positive measure. Put geometrically, it ensures that $S$ is not ``too flat" around any vertex.
    \item $k$ -- the maximum number of homological critical values for $f_{v}$ for any $v \in \mathbb{S}^{d-1}$, i.e. values at which the Euler characteristic of a sublevel set changes (assuming this quantity is finite). If one works with the PHT instead of the ECT, one considers homological critical values instead, where the homology of a sublevel set changes.
\end{itemize}

They show the following finiteness result:

\begin{theorem}[Thm. 7.1 in \cite{curry2018many}]
For either the ECT or the PHT, let $M \subset \mathbb{R}^d$ be a linearly  embedded simplicial complex, with appropriate constants $\delta, k$ as in the prior description. Then there is a constant $\Delta(d,\delta,k)$ such that $M$ can be determined using $\Delta(d,\delta,k)$ directions of the chosen transform.
\end{theorem}

The proof of this theorem is a multi-part algorithm, where the data computed at each step is passed forward as input to the next step. To begin, they show that, for a fixed $d$, an upper bound on $k$ and a lower bound on $\delta$ provide a bound on the total number of vertices in $M$ (Lemma 7.4 in \cite{curry2018many}). They then show that, given any sufficiently large collection of $\delta$-nets on the sphere, the resulting set of directions can be used to determine the location of the vertices in $M$ (Proposition 7.1 in \cite{curry2018many}). With the location of the vertices identified, one defines the following hyperplane arrangement in $\mathbb{R}^d$: $W(V) = \left[ \bigcup_{(v_1,v_2) \in \{V \times V - \Delta\}} (v_1 - v_2)^{T} \right]$, where $V$ is the vertex set of $M$, and where $\Delta$ is the diagonal in $V \times V$. That is, $W(V)$ is the union of all the hyperplanes in $\mathbb{R}^d$ orthogonal to the differences of pairs of distinct vertices in $V$. The connected components of $\mathbb{S}^{d-1} \cap \left( \mathbb{R}^d \setminus W(V) \right)$ are the $(d-1)$-dimensional strata of the stratification of the sphere induced by $W(V)$. The crucial observation to be made is that any two directions in the same top-dimensional stratum induce the same ordering on the simplices of $M$. Thus, given the ECT or PHT for any one direction in a stratum, it is possible to parametrize the ECT or PHT for all the other directions, provided the locations of the vertices $V$ are known (Lemma 5.3, Proposition 5.2 in \cite{curry2018many}). Thus, after identifying the set $V$ in the prior step, computing the hyperplane arrangement $W(V)$, and picking a test direction in each top-dimensional stratum, one has enough data to deduce the ECT or PHT on all of $\mathbb{S}^{d-1} \cap \left( \mathbb{R}^d \setminus W(V) \right)$, and, by continuity, on the entire sphere $\mathbb{S}^{d-1}$. Since the ECT or PHT on the full sphere determines the simplicial complex $M$ by prior injectivity results, we can ultimately deduce $M$ itself. The total number of directions needed in this procedure is
\[\Delta(d,\delta,k) = \left( (d-1)k \left( \frac{2\delta}{\sin(\delta)} \right)^{d-1} + 1 \right)\left( 1 + \frac{2}{\delta}\right)^{\delta} + O\left(\frac{dk}{\delta^{d-1}}\right)^{2d}\]

The proofs in both \cite{belton2018learning} and \cite{curry2018many} rely heavily on the simplicial complex structure of $M$, and there are presently no finiteness results known for more general subanalytic sets.

\subsection*{Sample ECT Code}

\begin{figure}[htb]
    \centering
    \includegraphics[scale=0.7]{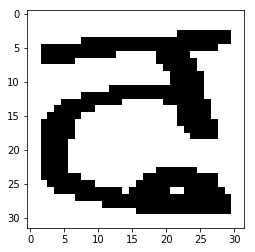}
     \includegraphics[scale=0.7]{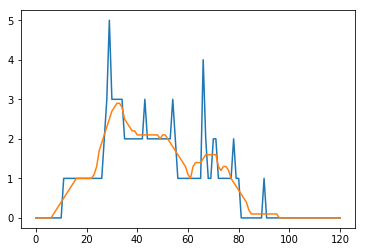}
    \caption{Left: greyscale image of a handwritten letter in the Devanagari alphabet, used in many North Indian languages. Right: the ECT of the above letter, taken in the direction $v = \langle 1, 1\rangle$. Superimposed on the ECT is a smoothed version, obtained via convolution. Images reproduced from \cite{Solomon2018}.}
    \label{fig:character}
\end{figure}

The author E. Solomon maintains a small GitHub repository with Python code for computing and comparing Euler Characteristic Transforms of 2D images \cite{Solomon2018}. The code samples the ECT along a finite set of directions for each image, and sets the distance between images to be the sum of the $L^2$ norms between smoothed Euler curves in matching directions. The choice and number of directions, smoothing parameter, and resulting classifier all have an impact on the prediction accuracy, although this is not well understood on a theoretical level at the moment. See Figure \ref{fig:character}.

\subsection{Intrinsic Persistent Homology Transform}

The PHT and ECT transforms of the prior section apply to shapes \emph{embedded} in $\mathbb{R}^d$. In \cite{BarcodeEmbeddings}, Oudot and Solomon propose an \emph{intrinsic} topological transform, the IPHT\footnote{This invariant is called the \emph{Barcode Transform} in that paper, but the name proposed above is clearer and fits better with existing literature}. This transform uses the extended persistence of a real-valued function $f: X \to \mathbb{R}$.

\begin{definition}
Let $(X,d_X)$ be a compact metric space. For each basepoint $p \in X$, consider the ``distance-to-the-basepoint" function $f_{p}(x) = d_{X}(p,x)$, and define $\Psi_{X}(p)$ to be the extended persistence of the pair $(X,f_p)$. We define $IPHT(X)$ to be the \emph{image} of $\Psi_{X}$, which, because $\Psi_{X}$ is continuous (a corollary of Theorem \ref{thm:funcstab}), is a compact subset of barcode space $\mathcal{B}$.
\end{definition}

It would seem that the appropriate analogue of the PHT of a subanalytic set $M$, as a map $PHT(M): \mathbb{S}^{d-1} \to \mathcal{B}$, would be the map $\Psi_{X}: X \to \mathcal{B}$. However, the map $\Psi_{X}$ has two shortcomings. Firstly, it requires us to keep track of the space $X$ as the domain of the map, when one would prefer a transform that allows us to discard the initial space $X$. Secondly, there is no simple way of comparing $\Psi_{X}$ and $\Psi_{Y}$ for distinct metric spaces $X$ and $Y$. This problem is resolved by taking $IPHT(X)$ to be the image of $\Psi_{X}$, so that it sits in a common ambient space for any choice of $X$.\\

So far, the IPHT has largely been studied in the context of \emph{compact metric graphs}, these being metric spaces arising from the shortest-path-metric on a weighted graph. In \cite{dey2015comparing}, Dey et al. propose the \emph{persistence distortion distance} $d_{PD}$ on the space of compact metric graphs $\mathbf{MGraphs}$.

\begin{definition}[\cite{dey2015comparing}]
 Let $d_{H}^{B}$ denote the Hausdorff distance on the space of compact subsets of the barcode space $\mathcal{B}$ induced by the Bottleneck distance. Then for any pair $X,Y \in \mathbf{MGraphs}$, define $d_{PD}(X,Y) = d_{H}^{B}(IPHT(X),IPHT(Y))$.
\end{definition}

\begin{figure}[htb]
	\includegraphics[scale = 0.4]{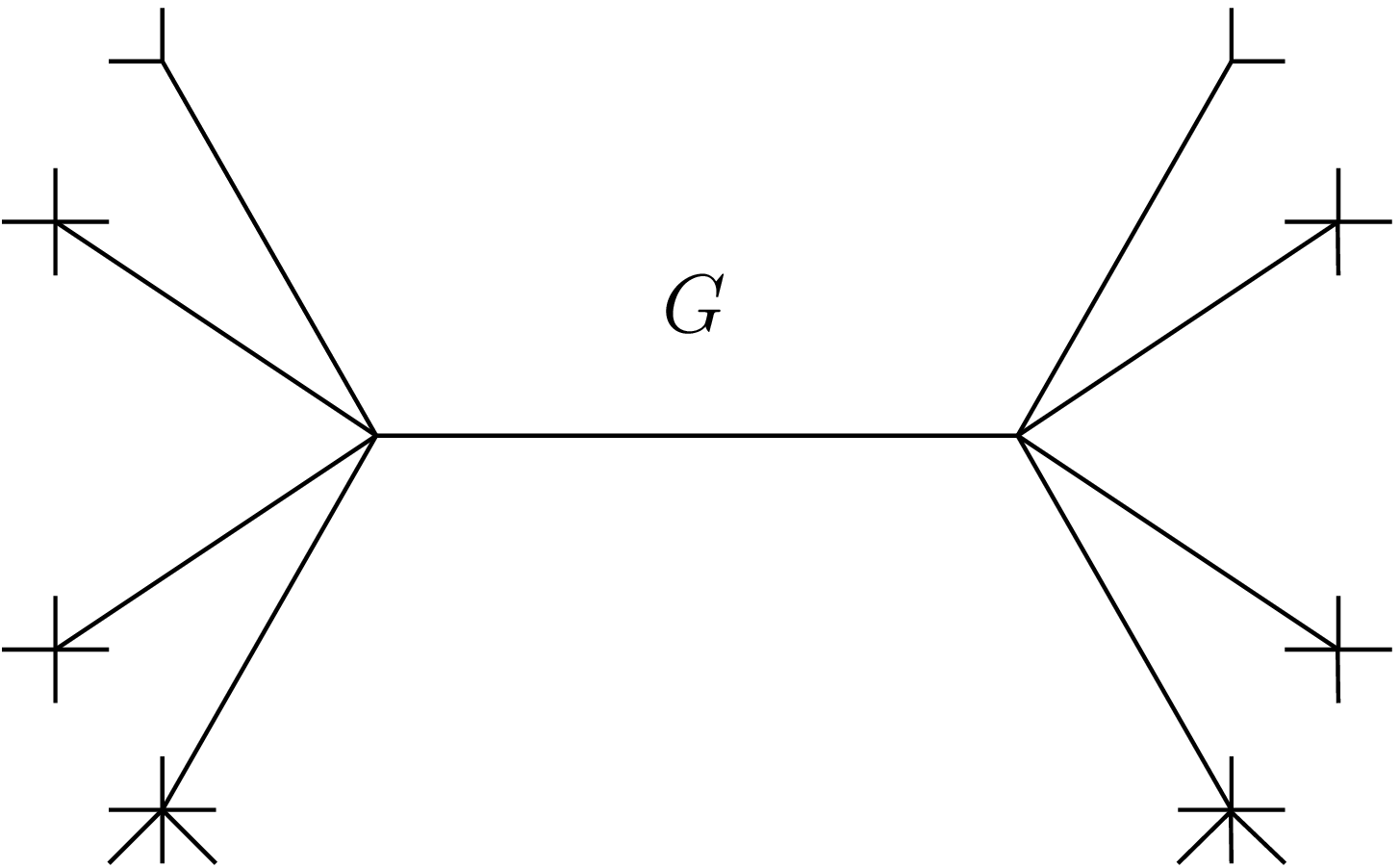} \,\,\,\,\,\,\,\,\,\,\,\,\,\,\,\,\,\,\,\,\,\,\,\,\,\,\,\,\,\,\,\,\,\,\,\,\,\,\,\,\,\,\,\,\,\,\,\,\,\,\,\,\,\,\,\,\,
	\includegraphics[scale = 0.4]{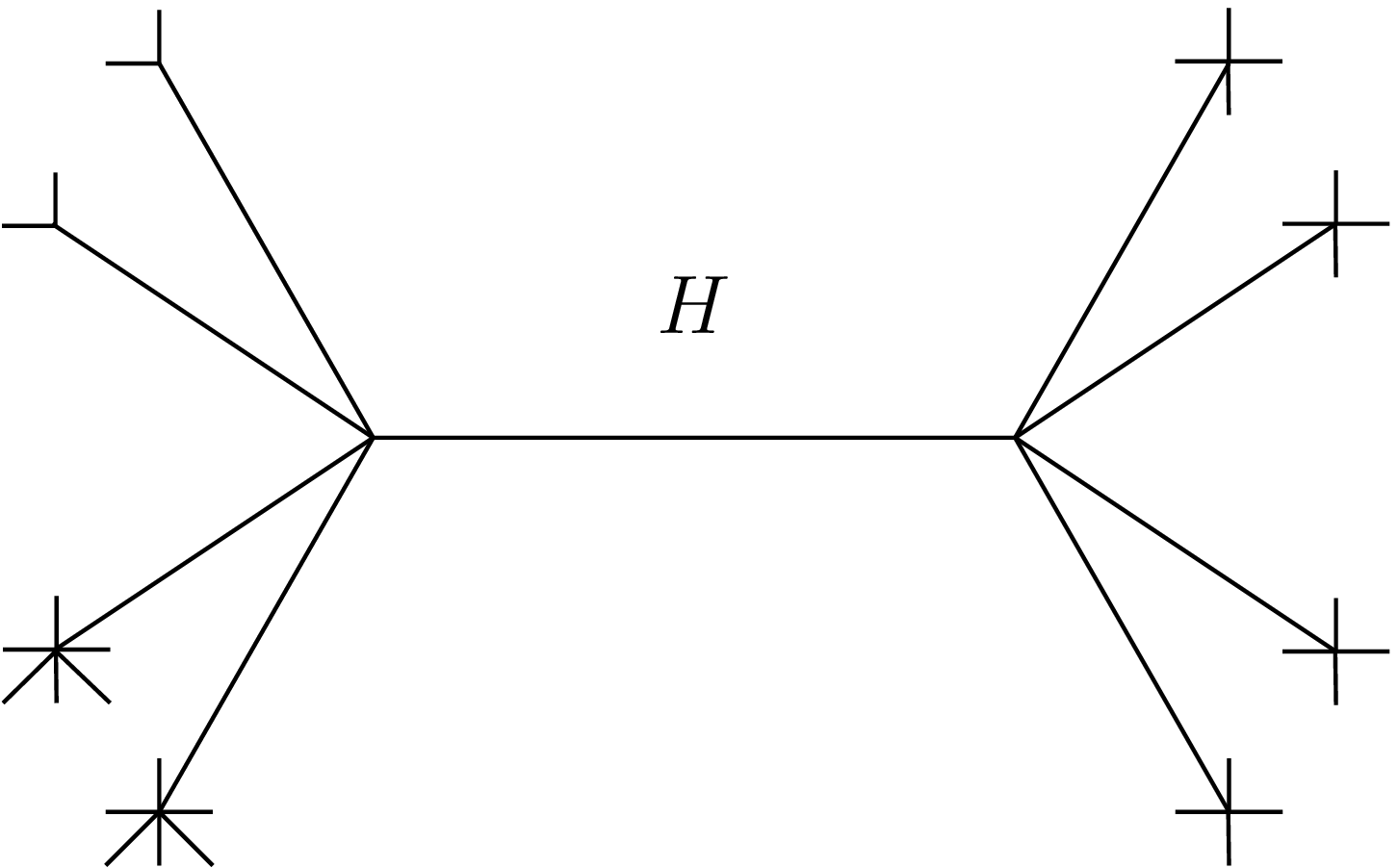}
	
	\caption{In the above figure, the lengths of the small branches are all equal to $1$, the lengths of the middle-sized branches are all equal to $10$, and finally both central edges have length $100$. For every middle-sized branch in $X$ there is a corresponding branch in $Y$ with the same number of small branches, not necessarily on the same side. The barcodes for points on matching branches are the same. Similarly the barcodes for points along the central edges of $X$ and $Y$ agree. Thus $\mathcal{BT}(X) = \mathcal{BT}(Y)$, but $X$ and $Y$ are not isomorphic.}
		\label{btnotinjfigure}	
\end{figure}

Dey et al.~\cite{dey2015comparing} show that the persistence distortion distance $d_{PD}$ is related to the Gromov-Hausdorff distance $d_{GH}$ as follows: $d_{PD}(X,Y) \leq 18 d_{GH}(X,Y)$ (see their Theorem 3 and the remark after their Theorem~4). In other words, the IPHT is a Lipschitz map on the space of metric graphs. In \cite{carriere2015local}, Carri\`{e}re et al. extend this result to a local Lipschitz property on the space of compact geodesic spaces. Dey et al. also show in~\cite{dey2015comparing} that $d_{PD}$ is computable in polynomial time (see their Theorem 26), provide an algorithm to do so, and conduct some preliminary experiments. Most relevant for our study of inverse problems, they demonstrate the existence of non-isometric graphs $X$ and $Y$ with $IPHT(X) = IPHT(Y)$. Oudot and Solomon~\cite{BarcodeEmbeddings} provide another (simpler) example of such a pair of graphs, see Figure \ref{btnotinjfigure} (coming from their Counterexample 5.2). This implies that the $IPHT$ is not an injective invariant on the space $\mathbf{MGraphs}$, and that $d_{PD}$ is only a pseudometric.\\

The pair of graphs $X$ and $Y$ in Figure \ref{btnotinjfigure}, as well as the pair of graphs in \cite{dey2015comparing}, have nontrivial automorphism groups. If $X$ has a nontrivial automorphism $\phi \in \operatorname{Aut}(X)$, then $\Psi_{X}(x) = \Psi_{X}(\phi(x))$ for all $x \in X$. This implies that the map $\Psi_{G}: G \to IPHT(G)$ is not injective, and hence we cannot recover the topological type of $X$ from that of $IPHT(X)$. In part, this can help explain the failure of injectivity of the $IPHT$. Now, Oudot and Solomon demonstrate that $\Psi_{X}$ can fail to be injective even if $\operatorname{Aut}(X)$ is trivial (Figure 6.1 in \cite{BarcodeEmbeddings}), so that injectivity of $\Psi_{X}$ is a stronger condition than $\operatorname{Aut}(X)$ being trivial. This motivated them to propose that the $IPHT$ might be injective on the set $\operatorname{INJ}_{\Psi}$ of graphs for which $\Psi_{X}$ is injective\footnote{Curry et al. make use of a similar genericity assumption when using the \emph{image} of the ECT as invariant: that every direction vector produces a distinct Euler curve~\cite[Def. 6.1]{curry2018many}}. They prove that this is indeed the case:

\begin{theorem}[Thm. 5.4 in \cite{BarcodeEmbeddings}]
The IPHT is injective up to isometry on the set $\operatorname{INJ}_{\Psi}$.
\label{thm:injfrominj}
\end{theorem}

The proof of Theorem \ref{thm:injfrominj} is based on the following pair of observations. On the one hand, Theorem \ref{thm:funcstab} implies that for any pair of points $x,x'$ in a metric graph $X$, $d_{B}(\Psi_{X}(x),\Psi_{X}(x')) \leq d_{X}(x,x')$. On the other hand, it is possible to show that for every $x \in X$ there exists a constant $\epsilon(x)$, such that if $x' \in G$ is another point with $d_{X}(x,x') \leq \epsilon(x)$, then $d_{B}(\Psi_{X}(x),\Psi_{X}(x')) \geq d_{G}(x,x')$. Taken together, these inequalities demonstrate that $\Psi_{X}$ is a \emph{local isometry}, i.e. for $d_{G}(x,x') \leq \epsilon(x)$, we have $d_{B}(\Psi_{X}(x),\Psi_{X}(x')) = d_{X}(x,x')$. Now, if $\Psi_{X}$ is injective, it is a homeomorphism onto $IPHT(G)$ (as its domain is compact and its codomain is Hausdorff). This implies that $IPHT(X)$ is homeomorphic to $X$. If we consider the \emph{intrinsic path metric} $\hat{d}_{B}$ on $IPHT(X)$ defined using the Bottleneck distance, the local isometry result then implies that $(IPHT(X),\hat{d}_{B})$ is globally isometric to $(X,d_{X})$. Thus, when $\Psi_{X}$ is injective, we have an explicit procedure for recovering $X$ from $IPHT(X)$, providing us with a left inverse\footnote{Note that when $\Psi_{G}$ is not injective, the set of continuous paths between points $x$ and $x'$ cannot be identified, via the map $\Psi_{X}$, with the set of continuous paths from $\Psi_{X}(x)$ to $\Psi_{X}(x')$. Thus, the local isometry result does not extend to a global isometry for the induced path metric.}.\\

The remainder of the paper demonstrates the extent to which the set $\operatorname{INJ}_{\Psi}$ is \emph{large} or \emph{generic}. For the Gromov-Hausdorff topology, they prove the following:

\begin{prop}[Prop. 5.5 in \cite{BarcodeEmbeddings}]
The set $\operatorname{INJ}_{\Psi}$ is Gromov-Hausdorff dense in $\mathbf{MGraphs}$.
\end{prop}

\begin{figure}[htb]
        \centering
        \includegraphics[scale = 0.3]{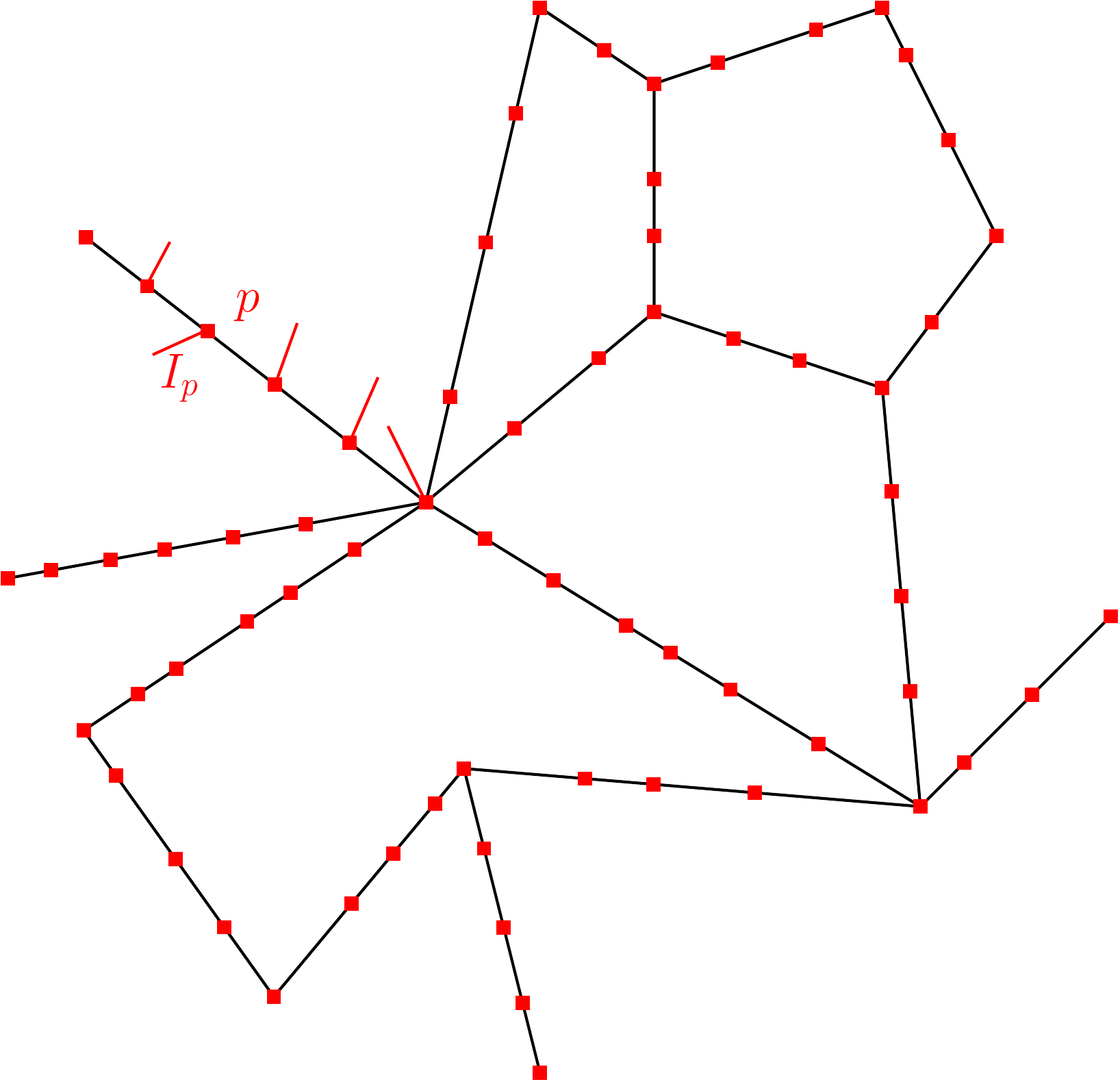}
        \caption{A graph $X$, drawn in black, with short \emph{thorns} of distinct lengths, drawn in red, attached along its edges. The resulting graph $X'$ is called a \emph{cactification} of $X$ in \cite{BarcodeEmbeddings}.}
        \label{fig:cactus}
    \end{figure}

The proof of this proposition is constructive: it demonstrates how to take a metric graph $X$ and insert many small branches along its edges so as to break any local geometric symmetry, forcing the map $\Psi$ to become injective (see Figure \ref{fig:cactus}).\\

In contrast, they show that every Gromov-Hausdorff open set admits a pair of non-isometric metric graphs $X$ and $X'$ with $IPHT(X) = IPHT(X')$ (Proposition 5.3 in \cite{BarcodeEmbeddings}), so that the IPHT cannot be injective on a Gromov-Hausdorff generic (open and dense) subset of $\mathbf{MGraphs}$. The proof of this result relies on finding an initial pair of non-isometric graphs with the same IPHT, shrinking them down, and gluing them to any other metric graph $Y$. In the Gromov-Hausdorff metric, the graphs $X$ and $X'$ that result from this gluing will be close to $Y$, and hence to each other, and it is not hard to show that $IPHT(X) = IPHT(X')$.\\

This suggests that the Gromov-Hausdorff topology is too coarse for studying this inverse problem, so the authors consider a finer topology on $\mathbf{MGraphs}$, the \emph{fibered topology}:

\begin{definition}[Def. 5.8 in \cite{BarcodeEmbeddings}]
For every combinatorial graph $X = (V,E)$, the set of metric structures on $X$ can be identified with the Euclidean fan $\mathbb{R}_{>0}^{E}/\operatorname{Aut}(X)$: one uses a vector in $\mathbb{R}_{>0}^{E}$ to assign edge weights, quotienting out by the automorphism group of $X$ to identify vectors of edge weights that produce isometric graphs. By restricting the focus to combinatorial graphs without valence-two vertices (which can be added to, or removed from, a graph without changing its topology), and by considering all possible combinatorial graphs satisfying this condition, one obtains a bijection $\gamma$ between $\mathbf{MGraphs}$ and the set $\Omega = \displaystyle\bigsqcup_{X = (V,E)}\mathbb{R}_{>0}^{E}/\operatorname{Aut}(X)$. Equipping each Euclidean fan with the topology induced on the quotient by the $L^2$ metric, one can then give $\Omega$ the disjoint union topology. Passing through the bijection $\gamma^{-1}$, one obtains a topology on $\mathbf{MGraphs}$ called the \emph{fibered topology}, as it decomposes that space into a countable family of disjoint open sets (the fibers).
\end{definition}

This fibered topology arises naturally when considering probability measures on $\mathbf{MGraphs}$ defined as mixture-models, where one first selects one of (countably many) combinatorial graphs $X = (V,E)$ and then chooses edge weights with a Borel measure on $\mathbb{R}_{>0}^{E}$ with density with respect to Lebesgue measure.\\

A lengthy combinatorial argument (Section 10 in \cite{BarcodeEmbeddings}) demonstrates that if $\Psi_{X}$ fails to be injective for some metric graph $X$, the set of edge lengths in $X$ is linearly dependent\footnote{To be precise, this only holds for graphs with at least three vertices, and for which there are no self-loops. A similar statement holds for the remaining cases, which is the focus of Section 11 in \cite{BarcodeEmbeddings}.} over $\mathbb{Z}$. Taking the contrapositive of this statement, one deduces that if the set of edge lengths in a graph $X$ is linearly independent over $\mathbb{Z}$, then $\Psi_{X}$ is injective. This linear independence condition is open and dense in the topology induced by the $L^2$ metric on each fiber $\mathbb{R}_{>0}^{E}/\operatorname{Aut}(X)$, and hence on all of $\mathbf{MGraphs}$ in the fibered topology. The authors thus conclude with the following injectivity result:

\begin{theorem}[Thm. 5.9A in \cite{BarcodeEmbeddings}]
There is a subset $U \subset \mathbf{MGraphs}$ containing $\operatorname{Inj}_{\Psi}$ on which the IPHT is injective, and which is generic in the fibered topology.
\end{theorem}

Oudot and Solomon also provide stability and injectivity results for a \emph{metric-measure} version of the IPHT (Theorems 4.2 and 5.9B in \cite{BarcodeEmbeddings}). They also prove the following Gromov-Hausdorff local injectivity result.

\begin{theorem}[Thm. 5.7 in \cite{BarcodeEmbeddings}]
	For every metric graph $X \in \mathbf{MGraphs}$ there is a constant $\epsilon(X) > 0$, such that if $Y$ is another compact metric graph with $0 < d_{GH}(X,Y) < \epsilon(X)$ then $d_{H}^{B}(IPHT(X), IPHT(Y)) >0$.
\end{theorem} 

\section{Conclusion}

The results explored in this survey form a preliminary but promising line of research into explainability from the topological point of view. Looking forward, there are a number of mathematical and data-theoretic challenges to be overcome:

\begin{itemize}
	\item How to best choose a set of direction vectors when implementing the PHT or ECT in practice.
	\item Finding the appropriate formulation of the IPHT that best extends to higher-dimensional intrinsic spaces, and investigating the associated injectivity properties (or lack thereof).
	\item Formulating the IPHT in terms of other families of functions defined on a metric space. For example, the barcodes arising from eigenfunctions of the Laplacian on metric simplicial complexes and manifolds.
	\item Studying the problem from the algorithmic point of view, including efficient implementations, bounds on complexity, etc.
	\item As of the writing of this survey, little is known about the statistics of the topological transforms and constructions detailed above. For example, the question of hypothesis testing has not been rigorously investigated.
\end{itemize}

Continued work by mathematicians, statisticians, and computer scientists will hopefully help address these questions, and bring the ideas discussed in this survey and their applications to maturity.

\phantomsection
\bibliographystyle{unsrt}
\bibliography{mybib}


\end{document}